\RequirePackage{xcolor}
\documentclass[a4paper]{article}

\usepackage[l2tabu, orthodox]{nag}
\usepackage{amsmath,amsfonts,amssymb,amsthm}
\usepackage{booktabs}
\usepackage[group-separator={,}]{siunitx}
\usepackage{fix-cm}
\usepackage{tikz}
\usepackage{pgfplots}
\usepackage[disable]{todonotes}
\usepackage{subcaption}
\usepackage{url}

\usepackage{fullpage}
\usepackage{lineno}
\usepackage[noabbrev]{cleveref}

\DeclareMathOperator*{\argmax}{arg\,max}

\pgfplotsset{compat=1.13}
\usetikzlibrary{external}
\usetikzlibrary{pgfplots.statistics}
\tikzsetexternalprefix{figures/}
\usepgfplotslibrary{statistics}
\graphicspath{{./img/}}

\makeatletter

\theoremstyle{definition}
\newtheorem{mydef}{Definition}[section]

\renewcommand{\todo}[2][]{\tikzexternaldisable\@todo[#1]{#2}\tikzexternalenable}

\def\input@path{{./img/},{./}}
\makeatother

\providecommand{\keywords}[1]{\noindent\textit{Keywords}: #1}
\providecommand{\subject}[1]{\noindent\textit{Subject areas}: #1}
\providecommand{\aucontribute}[1]{\paragraph*{Authors' contributions.}
  #1}
\providecommand{\funding}[1]{\paragraph*{Funding.} #1}
\providecommand{\competing}[1]{\paragraph*{Competing interests.} #1}

\newcommand*\samethanks[1][\value{footnote}]{\footnotemark[#1]}

\begin{document}

\title{A comparison of control strategies applied to a
  pricing problem in retail}

\author{A.~N.~Riseth\thanks{Mathematical Institute, University of
    Oxford, OX2 6GG, UK.},\hspace{0.1em}\thanks{Corresponding author, \texttt{riseth@maths.ox.ac.uk}.}
  \ J.~N.~Dewynne\samethanks[1], and C.~L.~Farmer\samethanks[1]}

%\subject{applied mathematics, revenue management}
%\keywords{dynamic pricing, stochastic control, approximate dynamic programming}

% \listoftodos

% \def\biblio{}
% \def\listoftodos{}

\maketitle

\begin{abstract}
  When sales of a product are affected by randomness in demand,
  retailers can use dynamic pricing strategies to maximise their
  profits. In this article the pricing problem is formulated as a
  stochastic optimal control problem, where the optimal policy can be
  found by solving the associated Bellman equation. The aim is to
  investigate Approximate Dynamic Programming algorithms for this
  problem. For realistic retail applications, modelling the problem
  and solving it to optimality is intractable. Thus practitioners make
  simplifying assumptions and design suboptimal policies, but a thorough
  investigation of the relative performance of these policies is
  lacking.

  To better understand such assumptions, we simulate the performance
  of two algorithms on a one-product system. It is found that for more
  than \emph{half of the realisations} of the random disturbance, the
  often-used, but approximate, Certainty Equivalent Control
  policy yields larger profits than an optimal, maximum expected-value
  policy. This approximate algorithm, however, performs significantly
  worse in the remaining realisations, which colloquially can be
  interpreted as a more risk-seeking
  attitude by the retailer. Another policy, Open-Loop Feedback
  Control, is shown to work well as a compromise between the Certainty
  Equivalent Control and the optimal policy.
\end{abstract}

\subject{Applied mathematics, Revenue management}

%%%%%%%%%%%%%%%%%%%%%%%%%%%
% Intro
%%%%%%%%%%%%%%%%%%%%%%%%%%

\section{Introduction}
The process of pricing products in order to control demand and
maximise revenues has been undertaken for centuries. In recent
decades, data- and model-driven approaches have become increasingly
popular in order to advise on and automate the process for companies.
There are several success stories from early adopters, for example in
the airline industry.
American Airlines estimated in 1992 that the introduction of revenue
management software had, over the preceding three
years, contributed \$500 million of additional revenue per year,
and would continue to do so in the future~\cite{smith1992yield}.
In an example from Chilean retail, the authors of~\cite{bitran1998coordinating}
report an expected revenue improvement of 7\%-16\% from implementing
model-driven strategies.
For retailers with billions of pounds in revenues, small
improvements to their revenue-management processes can be worth millions
of pounds.
In addition unsold items add up to thousands of tonnes of waste per year, so
improving the control of demand for products is advantageous
for both retailers and the environment.

We are interested in strategies to  set the prices of
products dynamically, and thus
formulate the problem in a stochastic optimal control framework.
Let $t$ denote equally spaced, discrete, dimensionless time points in $\{0,1,\dots,T\}$.
Product stock levels $\hat{S}={(\hat{S}_t)}_{t=0}^T$ are controlled by a pricing process
$\hat{\alpha}={(\hat{\alpha}_t)}_{t=0}^{T-1}$, and evolve according to a transition function $\hat{f}$ with random disturbance
$\hat{W}={(\hat{W}_t)}_{t=1}^T$,
\begin{equation}
  \hat{S}_{t+1}=\hat{f}(t,\hat{S}_t,\hat{\alpha}_t,\hat{W}_{t+1}).
\end{equation}
The goal is to find an $\hat{\alpha}$ which maximises the expected net revenues over a
given time horizon $T$,
\begin{equation}
  \max_{\hat{\alpha}}\mathbb E_{\hat{W}}\left[ \sum_{t=0}^T\hat{U}_t(\hat{S}_t,\hat{\alpha}_t,\hat{W}_{t+1})
    - \hat{\overline{U}}(\hat{S}_T)\right],
\end{equation}
where $\hat{U}_t$ denotes the revenue at time $t$, and $\hat{\overline{U}}$ the
cost of unsold stock.
We use the subscript on $\mathbb{E}_{\hat{W}}$ to emphasise that the
expectation is taken with respect to the random variables $\hat{W}_t$.
The solution to this problem can be found by solving the associated
Bellman equation, defined in~\eqref{eq:dynamic_programming_discrete}.

In practical applications in retail and other domains, it is almost always
intractable to solve the control problem to optimality~\cite{farmer2017uncertainty}.
The modelling of a real system and decision process can become very
complicated, and  we must  create a very
large state space in order to make use of the Bellman equation.
We often have to take into account unobservable
state variables, such as estimated parameters, and constraints
that may depend on history. Decision-makers in business may change their mind
about the objective over the course of the decision period, which should
be incorporated in the modelling as parameters with corresponding, estimated
probability distributions.
From a software implementation perspective, writing code that can
solve the Bellman equation efficiently can be much more complicated
than other methods.
Finally, the dimensions of the state, policy and exogenous
information spaces, can quickly make a numerical solution to the Bellman
equation intractable. An increase in the dimension can happen very
quickly, even for one-product problems.
In~\cite{bertsimas2001dynamic}, the authors model a simple
one-product demand function with uncertainty in the function
parameters, and present an eight-dimensional dynamic programming
solution to the problem.
Much of the focus in the research community has therefore been on
developing tractable algorithms with comparable average performance to the
optimal policies, see for
example~\cite{powell2011approximate} or~\cite{bertsekas2012dynamic}
for an overview.
One can either create explicit policy functions off-line, at the
start of the decision process, or implicitly through an automated
search for the best policy on-line for each decision point.
An advantage of the pre-calculated, explicit policy functions
is the speed at which we can make our decisions in the
future. In the classical pricing problems, the dynamics of the
underlying system do not require instant pricing decisions. Thus,
suboptimal decision rules that are created as they are needed, are often used
instead~\cite{talluri2006theory}. Estimation of the system and optimisation of prices are
normally separated, and constraints are more easily updated at each
decision point.

There are several proposed approximations in the
literature, although for revenue management applications, much of the domain knowledge is kept within
respective commercial organisations~\cite[Ch.~9]{talluri2006theory}.
Many of the suboptimal pricing algorithms are justified on practical
grounds~\cite{aviv2012dynamic}, or
from asymptotic results~\cite{gallego1994optimal}.
The aim of this article is to highlight the implications that suboptimal policy
choices have on the distribution of the relevant objective, and not
only marginalised quantities such as the expected value.
In particular, we consider a one-product pricing problem and
investigate the optimal Bellman policy, an often-used suboptimal
policy known as the \emph{Certainty Equivalent Control} (CEC)
policy, and a compromise between optimality and practicality known as the
\emph{Open-Loop Feedback Control} policy.
Rather surprisingly, we find that, more than half the
time, using the suboptimal CEC policy
leads to a higher profit than using the Bellman policy, for a wide
range of products. Of course, the
Bellman policy performs better on average, but
this is due to the CEC policy generating a larger, lower tail on the
distribution of profits than the Bellman policy.
Heuristically, one can say that the CEC policy has a higher risk
associated to it. We wish to emphasise that this indicates that a
choice of approximate control policy is implicitly a statement of
risk attitude. Thus, we propose that it should be addressed at the same level as a
description of a decision-maker's risk attitude as expressed with
utility functions or risk measures.

The article is organised as follows:
In \Cref{sec:bellman_optimal_control} we formulate the
pricing problem mathematically, and describe an algorithm to solve the
problem via the associated Bellman equation. We also give an example
of a problem with one product, and present the optimal pricing rules.
A discussion of the CEC and OLFC policies follows in
\Cref{sec:suboptimal_approximations}, where we compare its
performance to that of
the optimal policy.
Finally, we conclude and propose further work in \Cref{sec:conclusion}.

%%%%%%%%%%%%%%%%%%%%%%%%%%%
% Bellman
%%%%%%%%%%%%%%%%%%%%%%%%%%
\section{The optimal control problem}\label{sec:bellman_optimal_control}
We now consider a particular retail
pricing problem.
Given an initial amount of stock of a product and a future
termination time $T$,
the pricing problem is to
set the price
of the product dynamically, in order to both maximise the revenue and minimise the cost of
unsold stock at the termination time.
In this section, we define the optimal control pricing problem,
describe the optimality conditions given by the Bellman equation and
solve it numerically for an example system.

Consider a system over discrete, equispaced, dimensionless time points
$t=0,1,\dots,T$, with state
$\hat{S}_t$ and pricing process $\hat{\alpha}$ such that the $\hat{\alpha}_t$ that take values in a closed
interval $\hat{A}=\nobreak[\hat{a}_{\mathrm{min}},\hat{a}_{\mathrm{max}}]$ of prices.\footnote{In
  practice, policy values may have additional constraints that depend
  on the current state,
  in which case we say that $\hat{\alpha}_t$ must take values
  in a set $\hat{A}(\hat{S}_t)$.
}
The state $\hat{S}_t$ is the stock of a product at time $t$, and
$\hat{\alpha}_t$ the price for the product in the time period from $t$ to
$t+1$. We model the amount of product sold over each time period according
to a forecast demand
function $\hat{q}:\hat{A}\to\mathbb R_+$, which is bounded, continuous and decreasing.
Exogenous influences on demand are considered as randomness in the
system, and are
modelled in a multiplicative and dimensionless fashion by a stochastic process
$W=(W_1,\dots,W_T)$ taking non-negative values. See, for
example,~\cite[Ch.~7]{talluri2006theory} for a discussion of
demand models and the modelling of uncertainty.

For a given pricing process $\hat{\alpha}$, the
system evolves from some initial state $\hat{S}_0>0$, according to the
recursion
\begin{equation}\label{eq:stock_dynamics}
  \hat{S}_{t+1}^{\hat{\alpha}}=\hat{S}_t^{\hat{\alpha}}-\min(\hat{S}_t^{\hat{\alpha}},\hat{q}(\hat{\alpha}_t)W_{t+1}),\qquad t=0,\dots,T-1.
\end{equation}
The function $\hat{Q}(\hat{s},\hat{a},w)=\min(\hat{s},\hat{q}(\hat{a})w)$ denotes the unit sales over a
period at price $\hat{a}$,
starting with stock $\hat{s}$, and with exogenous influences
characterised by $w$. The minimum operator is used ensure that the
amount of product sold over the time period does not exceed the
current stock level.

The revenue accrued over period $t\to t+1$ is $\hat{\alpha}_t\hat{Q}(\hat{S}_t^{\hat{\alpha}},\hat{\alpha}_t,W_{t+1})$.
The cost of remaining stock at time $T$ is modelled by a cost per unit
stock $\hat{C}\geq 0$.\footnote{In some situations, unsold items at time $T$
  may be sold at some ``salvage price'', in which case we could allow
  $\hat{C}<0$. We assume $\hat{C}\geq 0$ in this paper.}
Let $\hat{\mathcal{A}}$ denote the set of feasible processes $\hat{\alpha}$ that take values
in $\hat{A}$.
Define the value of having stock $\hat{s}$ at time $t\leq T$
by the value function $\hat{v}$, such that
\begin{align}\label{eq:value_function_def}
  \hat{v}(t,\hat{s})&=\max_{\hat{\alpha}\in\hat{\mathcal{A}}} \hat{J}(t,\hat{s},\hat{\alpha}),\quad\text{where}\\
  \hat{J}(t,\hat{s},\hat{\alpha})&=
                                   \mathbb E_{W}\left[ \sum_{\tau=t}^{T-1}
                                   \hat{\alpha}_\tau \hat{Q}(\hat{S}_\tau^{\hat{\alpha}},\hat{\alpha}_\tau,W_{\tau+1})
                                   - \hat{C}\hat{S}_T^{\hat{\alpha}} \mid \hat{S}_t^{\hat{\alpha}} = \hat{s}
                                   \right].
                                   \label{eq:value_function_def2}
\end{align}
This leads us to the following mathematical formulation of the pricing
problem:
\begin{mydef}
  Given an initial amount of stock $\hat{S}_0>0$ and a cost per unit unsold stock $\hat{C}\geq
  0$, the \emph{pricing problem} is to find a pricing process $\hat{\alpha}^*\in\hat{\mathcal{A}}$ such that
  \begin{equation}
    \hat{J}(0,\hat{S}_0,\hat{\alpha}^*) = \hat{v}(0,\hat{S}_0),
  \end{equation}
  that is, $\hat{\alpha}^* = \argmax_{\hat{\alpha}\in\hat{\mathcal{A}}}\hat{J}(0,\hat{S}_0,\hat{\alpha})$.
\end{mydef}
We choose to maximise the expected profit over the period, which
assumes a risk-neutral decision-maker.\footnote{For characterisations of
  decision-makers and their risk attitude, see for
  example~\cite[App.~G]{bertsekas2005dynamic}.}
It is still
important, however,
to understand the distribution of profits for a given pricing policy
$\hat{\alpha}$. Therefore, we simulate the distribution when
we investigate the performance of algorithms in
\Cref{sec:suboptimal_approximations}. %,sec:markdown_miss_specification}.
If we assume that the random variables $W_t$ are independent, this
stochastic optimal control problem can be solved by
considering the optimality conditions that arise from the Dynamic
Programming principle, also known as the Bellman equation.

\subsection{Non-dimensionalisation of the system}
We will now consider a non-dimensionalised
representation of the system.
% The unit of time is already scaled such that
% the difference between two time points is $1$.
The units of stock will be  scaled with
respect to
the initial stock $\hat{S}_0$, and units of money will be scaled with respect to
the upper bound on
price, $\hat{a}_{\mathrm{max}}$.
Let the corresponding dimensionless quantities be defined without hats.
Then we set
\begin{align}
  s
  &= \frac{\hat{s}}{\hat{S}_0},
  & a
  &=\frac{\hat{a}}{\hat{a}_{\mathrm{max}}},
  &C&=\frac{\hat{C}}{\hat{a}_{\mathrm{max}}}.%,
  % &\hat t &= t,
  % &\hat W_{\hat t}&=W_{t}.
\end{align}
% To be precise, the dimensional time $t$ should be divided by $1$ unit of time.
The dimensionless functions $q(a)$ and $Q(s,a,w)$,
for forecasted demand
and realised sales respectively, are
\begin{align}
  q(a)&= \frac{\hat{q}(a\cdot \hat{a}_{\max})}{\hat{S}_0},
  &Q(s,a, w)&= \min(s, q(a) w).
\end{align}
The collection of pricing policies $\mathcal{A}$ contains all
the processes $\hat{\alpha}_t/\hat{a}_{\mathrm{max}}$, where
$\hat{\alpha}\in \hat{\mathcal{A}}$.
Now we can define the dimensionless value function
\begin{align}\label{eq:value_function_def_nondim}
  v(t,s)&=\max_{\alpha\in\mathcal{A}}
          J(t,s,\alpha),\quad\text{where}\\
  J(t,s,\alpha)&=
                 \mathbb{E}_{W}\left[ \sum_{\tau=t}^{T-1}
                 \alpha_\tau Q(
                 {S}_\tau^\alpha,\alpha_\tau,
                 W_{\tau+1})
                 - C S_{T}^\alpha \big\vert S_{t}^\alpha =
                 s
                 \right].
                 \label{eq:value_function_def_nondim2}
\end{align}
Finally, the dimensionless optimal control problem
is to find $\alpha^*\in{\mathcal A}$, such that
$J(0,1,\alpha^*)=v(0,1)$.
For the remainder of this article, we work with the
non-dimensionalised system.%, and drop the hats from the
%dimensionless quantities.

\subsection{The Bellman equation}
We choose $\mathcal A$ to be the set of Markovian policies for the
problem. This means that for each $\alpha\in \mathcal
A$, there exists a
measurable function $a:[0,T)\times \mathbb R_+\to A$, such that for each possible outcome
$\omega$, the process $\alpha$ is given by
$\alpha_t(\omega) =
a(t,S_t^\alpha(\omega))$.
Then,
an approach to finding the value function above is to use the Dynamic
Programming principle~\cite{bertsekas2005dynamic}, which states that $v$ can be defined
recursively by
\begin{equation}\label{eq:dynamic_programming_discrete}
  v(t,s)=\max_{a\in A}\mathbb E_{W}\left[
    aQ(s,a,W_{t+1})
    +v\bigl(t+1,s-Q(s,a,W_{t+1})\bigr)\right].
\end{equation}
Thus, the value function is the solution to the backwards-in-time
recursive relation~\eqref{eq:dynamic_programming_discrete} with
terminal value $v(T,s)=-Cs$, and the optimal policy
function $a(t,s)$ is given by the argmax for each $(t,s)$.
The recursion~\eqref{eq:dynamic_programming_discrete} is called the
\emph{Bellman equation}, and a discussion
of its validity can be found, for example, in~\cite{bertsekas2005dynamic}.
Analytical solutions to Bellman equations are only available in very
rare cases, and  numerical approaches are normally needed to
approximate the solutions.

We implement the following algorithm to solve the optimal control
problem, using the Bellman equation:\todo{Decorate with some algorithm environment?}
\begin{enumerate}
\item Create a grid of equispaced points $0=s_1<s_2<\cdots<s_K=S_0$, and arrays $v^K\in\mathbb R^{K\times(T+1)}$,
  $\alpha^K\in\mathbb R^{K\times T}$.
\item Set $I[v^K](s)=-Cs$.
\item Set $v^K[i,T]=I[v^K](s_i)$ for $i=1,\dots, K$.
\item For $t = T-1,\dots,0$:
  \begin{enumerate}
  \item Set $\displaystyle v^K[i,t]=\max_{a\in A}\mathbb E_{W_{t+1}}\left[ aQ(s_i,a,W_{t+1})
      +I[v^K](s_i-Q(s,a,W_{t+1}))\right]$\\ for $i=1,\dots,K$.
  \item Set $\alpha^K[i,t]$ to the maximiser above.
  \item Set $I[v^K](s) = \mathrm{Interpolate}(s, {(s_i)}_i,{(v^K[i,t])}_i)$.
  \end{enumerate}
\item Return $v^K,\alpha^K$.
\end{enumerate}
The expectation above is approximated using Monte--Carlo simulation,
for which the variance of the approximation error decreases by one over the
number of samples used~\cite{caflisch1998monte}.
In this article we use \num{1000} samples for the approximation of the
expectation in step~(iv)(a).
We choose to use piecewise
linear interpolation for $I[v^K](s)$ in step~(iv)(c).

\subsection{Example system}\label{sec:bellman_example_markdown}
\begin{figure}[htbp]
  \centering
  \begin{subfigure}[b]{0.5\textwidth}
    \begin{tikzpicture}[scale=0.8]
      \begin{axis}[
        xlabel={$s$},
        ylabel={$v(t,s)$},
        title={Value function},
        legend cell align=left,
        legend pos=north west
        ]
        \addplot+[thick,mark=none] table[x index = 0,y index = 1,col sep=comma]
        {./data/markdown_bellman_det_val_policy.csv};
        \addlegendentry{$t=0$};
        \addplot+[thick,mark=none,dashed] table[x index = 0,y index = 2,col sep=comma]
        {./data/markdown_bellman_det_val_policy.csv};
        \addlegendentry{$t=1$};
        \addplot+[thick,mark=none,dash dot] table[x index = 0,y index = 3,col sep=comma]
        {./data/markdown_bellman_det_val_policy.csv};
        \addlegendentry{$t=2$};
      \end{axis}
    \end{tikzpicture}
  \end{subfigure}%
  \begin{subfigure}[b]{0.5\textwidth}
    \begin{tikzpicture}[scale=0.8]
      \begin{axis}[
        xlabel={$s$},
        ylabel={$a(t,s)$},
        title={Policy function},
        legend cell align=left,
        ]
        \addplot+[thick,mark=none] table[x index = 0,y index = 5,col sep=comma]
        {./data/markdown_bellman_det_val_policy.csv};
        \addlegendentry{$t=0$};
        \addplot+[thick,mark=none,dashed] table[x index = 0,y index = 6,col sep=comma]
        {./data/markdown_bellman_det_val_policy.csv};
        \addlegendentry{$t=1$};
        \addplot+[thick,mark=none,dash dot] table[x index = 0,y index = 7,col sep=comma]
        {./data/markdown_bellman_det_val_policy.csv};
        \addlegendentry{$t=2$};
      \end{axis}
    \end{tikzpicture}
  \end{subfigure}
  \caption{The value function and corresponding
    optimal control function for the pricing problem.
    The regions where the value function moves from a linear to a
    nonlinear regime
    corresponds to when the pricing policy moves from the upper bound
    \num{1} to less than \num{1}.
  }\label{fig:markdown_bellman}
\end{figure}
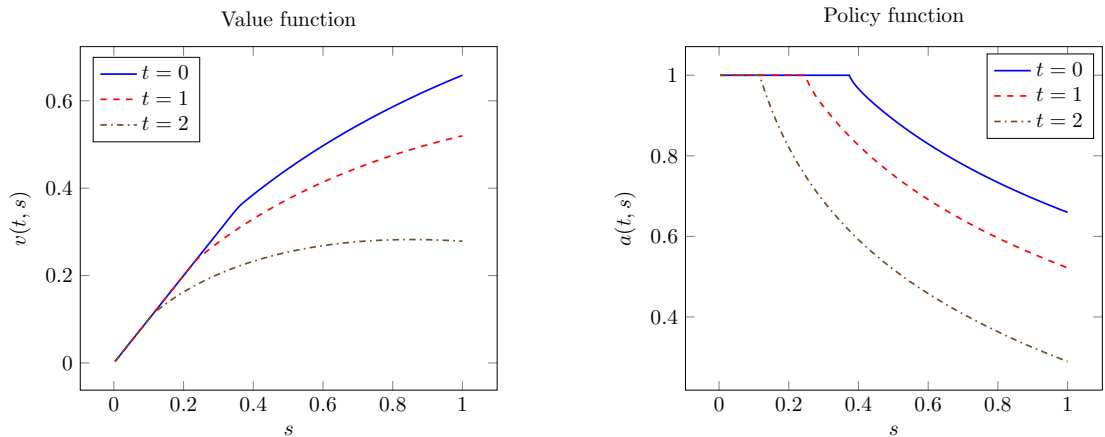

\begin{figure}[htbp]
  \centering
  \begin{subfigure}[t]{0.5\textwidth}
    \begin{tikzpicture}[scale=0.8]
      \begin{axis}[
        ylabel={Price},
        xlabel={\phantom{$P(\alpha)$}},% \phantom: to align the two figures
        title={$\alpha_t$, simulated values},
        boxplot/draw direction=y,
        xtick={1,2,3},
        xticklabels={$t=0$, $t=1$, $t=2$},
        ]
        \addplot+[boxplot] table[y index=0,col sep=comma]
        {./data/markdown_bellman_det_policies.csv};
        \addplot+[boxplot = {whisker range = 10}] table[y index=1,col sep=comma]
        {./data/markdown_bellman_det_policies.csv};
        \addplot+[boxplot = {whisker range = 10}] table[y index=2,col sep=comma]
        {./data/markdown_bellman_det_policies.csv};
      \end{axis}
    \end{tikzpicture}
  \end{subfigure}%
  \begin{subfigure}[t]{0.5\textwidth}
    \begin{tikzpicture}[scale=0.8]
      \begin{axis}[
        xlabel=Profit $P(\alpha)$,
        ylabel=Count,
        title={Realised profit},
        legend cell align=left,
        xmin=0.565, xmax=0.69,
        ymax=1400,
        ]
        \addplot[blue,hist={data=x,bins=40}] table [y index = 0, col
        sep=comma]
        {./data/markdown_bellman_det_vals.csv};
      \end{axis}
    \end{tikzpicture}
  \end{subfigure}%
  \caption{Simulations of the pricing system, started at $S_0=1$
    and controlled by the optimal policy $\alpha$ shown in
    \Cref{fig:markdown_bellman}.
    The left figure is a box plot
    that shows the median, the quartiles and the extremal prices observed
    in the simulations.
    As we see, the variance of the prices increases
    in time, reflecting the wider range of realised remaining stock at these times.
    The right hand figure is a histogram of the total profit over
    the pricing period.
  }\label{fig:bellman_simulation}
\end{figure}
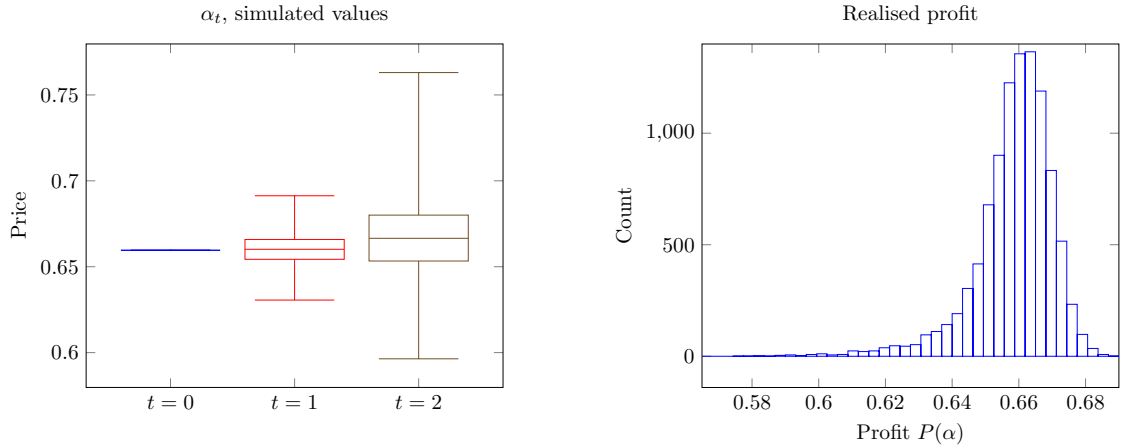

In order to investigate the optimal pricing of a specific system, we
choose to look at a family of demand functions of the form
$q(a)=q_1e^{-q_2a}$, where $q_1,q_2>0$. For a discussion
about their properties and usage in modelling demand, see~\cite[Ch.~7]{talluri2006theory}.
All the numerical experiments in this article use this family of exponential
demand functions, with price constraints $A=[0,1]$.
In the current section, and in
\Cref{fig:markdown_bellman,fig:bellman_simulation,fig:bellman_det_policy_difference,fig:bellman_det_vals},
we consider a particular choice of $q_1,q_2$ so that the demand
function is given by $q(a)=\frac{1}{3}e^{2-3a}$.
In \Cref{fig:profit_diff_heatmaps} and \Cref{tbl:paramcomparisons} a larger
range of values for $q_1,q_2$ are considered.

We assume the exogenous disturbance process is a sequence of
shifted, independent and identically Beta-distributed random
variables with mean $1$ and variance
$\gamma^2$. That is, we define $W_t\sim \frac{1}{2}+X$, where
$X\sim \mathrm{Beta}(\mu,\nu)$,\footnote{A
  $\mathrm{Beta}(\mu,\nu)$ random variable has probability density
  function on $[0,1]$, given by
  $p(x)=\frac{x^{\mu-1}{(1-x)}^{\nu-1}}{B(\mu,\nu)}$.
  The function $B(\mu,\nu)$ is the normalising factor, and we require
  that $\mu,\nu>1$.
}
and $\mu=\nu=\frac{1}{8\gamma^2}-\frac{1}{2}$.
To ensure that $X$ is unimodal, we require that $\gamma^2<1/12$.

Set $\gamma = 5\times 10^{-2}$, $C=1$ and $T=3$.
The numerical solution to the Bellman equation, and the corresponding optimal
pricing policy, is shown in \Cref{fig:markdown_bellman}.
Let us now investigate the behaviour of the optimal pricing policy
$\alpha$ and
the outcome of following this policy.
Define a random variable $P(\alpha)$, which for each realisation
represents the total profit, is
\begin{equation}
  P(\alpha) = \sum_{t=0}^{T-1}\alpha_tQ(S_t^\alpha,\alpha_t,W_{t+1}) - CS_T^\alpha.
\end{equation}
By sampling from the stochastic process $W=(W_1,\dots,W_T)$, we can
estimate the random variables $\alpha_t$ and $P(\alpha)$.
The plots in \Cref{fig:bellman_simulation} show the results of
simulating the system \num{10000} times.\footnote{Preliminary
  investigations indicated that \num{1000}
  simulations is sufficient to obtain an accurate representation of
  the distributions in this article. However, we have chosen to use
  more when this is computationally convenient.
}

%%%%%%%%%%%%%%%%%%
% Comparison
%%%%%%%%%%%%%%%%%
\section{Suboptimal
  approximations}\label{sec:suboptimal_approximations}
In this section we will look at two suboptimal policies among the class of
algorithms that calculate decisions on-line.
Whenever a decision must be made, these methods simulate the future
behaviour of the system and optimise the current decision based on
these simulations.
We consider two special cases known as the
\emph{Certainty Equivalent Control} policy, which uses a point
estimate of the system, and
the \emph{Open-Loop Feedback Control} policy, which uses more
information about the  future behaviour of the system~\cite[Ch.~6]{bertsekas2005dynamic}.
For the example system in this article, it turns out that we can find
the CEC policy analytically.
Numerical comparison experiments between these policies and the optimal
Bellman policy are
made using the example pricing problem introduced in
Subsection~\ref{sec:bellman_optimal_control}\ref{sec:bellman_example_markdown}. We also extend the comparison to
cover a much larger range of parameters in the problem, to show that
the same results hold for a wide range of retail products.

\subsection{The Certainty Equivalent Control policy}
An often-used, tractable, algorithm for solving stochastic optimal control
problems is the Certainty Equivalent Control (CEC) policy, as it is
named in~\cite{bertsekas2005dynamic}.
It is also known as (deterministic) Model
Predictive Control in the engineering community.
The algorithm is particularly practical, because the optimisation
problem is reduced to a standard deterministic problem that can
be solved with existing, commercial or open-source solvers which
handle very large systems.
At each decision time, a deterministic optimisation problem for the
remaining decision horizon is solved. Only the decision for the
current time step is
used, whilst the subsequent decisions are discarded.
For each decision time $t=0,1\dots,T-1$, the CEC policy for the pricing
problem calculates the current price using the following steps:
\begin{enumerate}
\item Observe the state $s$.
\item Take a best estimate $w_{t+1},\dots,w_T$ of ${(W_\tau)}_{t+1}^T$.
\item Solve the optimisation problem
  \begin{equation}\label{eq:cec_optim_problem}
    \max_{\mathbf a\in A^{T-t}}\left\{\sum_{\tau=t}^{T-1}\mathbf
      a_\tau Q(S_\tau^{\mathbf a},\mathbf
      a_\tau,w_{\tau+1})-CS_T^{\mathbf a}\right\},
    \quad \text{s.t.}\quad S_t^{\mathbf a}=s,
  \end{equation}
  where $A^{T-t}$ is the $T-t$ times Cartesian product $A\times
  \cdots\times A$.
\item Set the price corresponding to the element
  $\mathbf{a}_t\in A$ of a maximiser to~\eqref{eq:cec_optim_problem}.
\end{enumerate}

For the one-product pricing problem, we can simplify the optimisation
problem and in some cases obtain analytical solutions for the policy function.
First, let us consider the expected value estimate $w_\tau=\mathbb E
[W_\tau]=1$ for $\tau=t+1,\dots T$.
Then we can rewrite the maximisation problem to find
the certainty equivalent value function,
\begin{equation}
  \widetilde{v}(t,s)=
  \max_{\mathbf a\in A^{T-t}}\left\{\sum_{\tau=t}^{T-1}(\mathbf
    a_\tau+C)\min\left(s-\sum_{r=t}^{\tau-1}q(\mathbf a_r),q(\mathbf a_\tau)\right)-Cs\right\}.
  % \quad \text{s.t.}\quad \sum_{\tau=t}^{T-1}q(\mathbf a_t)\leq s.
\end{equation}
The optimal choice here is to let $\mathbf a_\tau=a^*\in A$ for each
$\tau$, such that the same amount of stock is forecast to be sold in
each period.
Then, $a^*=a^C(t,s)$ is given by the policy function
\begin{equation}
  a^C(t,s)=\argmax_{a\in A} \left\{(a+C)\min\left(
      \frac{s}{T-t},q(a)
    \right)\right\}.
  % ,\quad\text{s.t.}\quad q(a)\leq \frac{s}{T-t}.
\end{equation}
Let $\mathcal P_A$ be the projection operator from $\mathbb R$ onto the interval $A$.
Let us consider two families of demand functions that are
popular in the literature, see for example~\cite[Ch.~7]{talluri2006theory}.
% If $q$ is of the form $a\mapsto q_1-q_2a$, with $q_1,q_2> 0$, then
% the optimal policy $\alpha_t=a^C(t,S_t^\alpha)$ for the deterministic
% problem is given by the function
% \begin{equation}
%   a^C(t,s)=\mathcal P_A \left[ \max\left(
%       \frac{q_1}{q_2}-\frac{s}{q_2(T-t)},\frac{1}{2}\left(\frac{q_1}{q_2}-C
%       \right) \right) \right].
% \end{equation}
In the case when $q(a)=q_1e^{-q_2a}$, the policy function is
\begin{equation}\label{eq:cec_policy}
  a^C(t,s)=\mathcal P_A\left[
    \max\left( \frac{1}{q_2}\log\left( \frac{q_1(T-t)}{s}\right),
      \frac{1}{q_2}-C  \right)\right].
\end{equation}
We start with investigating policy functions for a
single combination of the  system
parameters, and then provide a
comparison between the Bellman policy and CEC policy for a larger
range of parameters in Subsection~\ref{sec:suboptimal_approximations}\ref{sec:parameter_comparison}.

\subsubsection{Example system}\label{sec:cec_comparison_example}
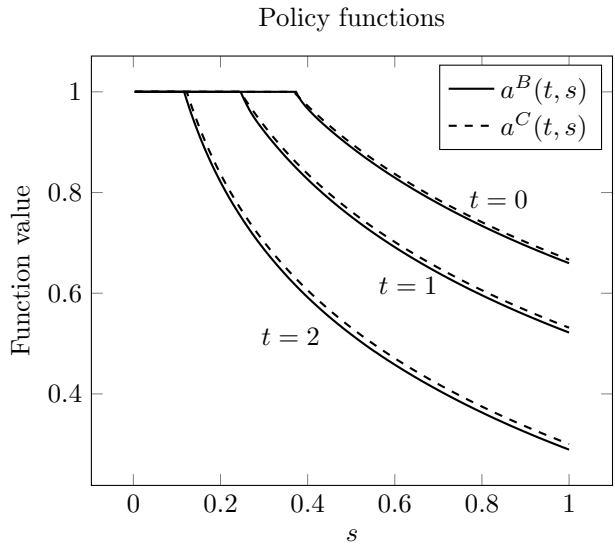
\begin{figure}[htbp]
  \centering
  \begin{tikzpicture}
    \begin{axis}[
      xlabel={$s$},
      ylabel={Function value},
      title={Policy functions},
      legend cell align=left,
      ]
      \addplot[thick,mark=none] table[x index = 0,y expr={\thisrowno{5}},col sep=comma]
      {./data/markdown_bellman_det_val_policy.csv};
      \addlegendentry{$a^B(t,s)$};
      \addplot[thick,mark=none,dashed] table[x index = 0,y expr={\thisrowno{8}},col sep=comma]
      {./data/markdown_bellman_det_val_policy.csv};
      \addlegendentry{$a^C(t,s)$};
      \addplot[thick,mark=none] table[x index = 0,y expr={\thisrowno{6}},col sep=comma]
      {./data/markdown_bellman_det_val_policy.csv};
      \addplot[thick,mark=none,dashed] table[x index = 0,y expr={\thisrowno{9}},col sep=comma]
      {./data/markdown_bellman_det_val_policy.csv};
      \addplot[thick,mark=none] table[x index = 0,y expr={\thisrowno{7}},col sep=comma]
      {./data/markdown_bellman_det_val_policy.csv};
      \addplot[thick,mark=none,dashed] table[x index = 0,y expr={\thisrowno{10}},col sep=comma]
      {./data/markdown_bellman_det_val_policy.csv};
      \draw (axis cs:0.45,0.55) node[anchor=north east] {$t=2$};
      \draw (axis cs:0.72,0.65) node[anchor=north east] {$t=1$};
      \draw (axis cs:0.75,0.75) node[anchor=south west] {$t=0$};
    \end{axis}
  \end{tikzpicture}

  % \vspace{1em}
  % \begin{tikzpicture}
  %   \begin{axis}[
  %     xlabel={$s$},
  %     ylabel={$[a^B(t,s)-a^C(t,s)]/a^B(t,s)$},
  %     title={\hspace{1em}Relative function difference},
  %     legend cell align=left,
  %     ]
  %     \addplot+[mark=none] table[x index = 0,y expr={1-\thisrowno{8}/\thisrowno{5}},col sep=comma]
  %     {./data/markdown_bellman_det_val_policy.csv};
  %     \addlegendentry{$t=0$};
  %     \addplot+[mark=none] table[x index = 0,y expr={1-\thisrowno{9}/\thisrowno{6}},col sep=comma]
  %     {./data/markdown_bellman_det_val_policy.csv};
  %     \addlegendentry{$t=1$};
  %     \addplot+[mark=none] table[x index = 0,y expr={1-\thisrowno{10}/\thisrowno{7}},col sep=comma]
  %     {./data/markdown_bellman_det_val_policy.csv};
  %     \addlegendentry{$t=2$};
  %   \end{axis}
  % \end{tikzpicture}
  % \todo[inline]{Remove the relative difference comparison?}
  \caption{Comparison of the CEC and Bellman policy functions $a^C$ and
    $a^B$.
    The CEC function sets
    higher prices than the Bellman function for large $s$.
    % Observe the spikes in the difference plot for $t<T-1$:
    % In a small interval, $a^B(t,s)\geq a^C(t,s)$, before both functions
    % reach the boundary price $a_{\mathrm{max}}=1$. This sudden, relative increase in
    % the Bellman function compared to $a^C$ is due to the
    % diffusion of the value function caused by the
    % random variables $W_t$.
  }\label{fig:bellman_det_policy_difference}
\end{figure}

The example system in Subsection~\ref{sec:bellman_optimal_control}\ref{sec:bellman_example_markdown} has a demand
function of the form $q(a)=q_1e^{-q_2a}$, where $q_1=e^2/3$ and
$q_2=3$. We can therefore compare the optimal Bellman policy with the
pricing policy that a CEC algorithm would imply.
Denote the Bellman policy function by $a^B$, and consider
$a^B(t,s)-a^C(t,s)$ for each $t=0,\dots,T-1$. The plot in
\Cref{fig:bellman_det_policy_difference} shows this difference.
Both of the policy functions reach the upper bound in $A$ for small values
of $s$, but most of the time, the Bellman policy is pricing the products
lower than the CEC policy.

\begin{figure}[htbp]
  \centering
  \begin{tikzpicture}[scale=0.8]
    \begin{axis}[
      xlabel=Profit $P(\alpha^B)$,
      ylabel=Count,
      title={Simulations of Bellman policy},
      xmin=0.565, xmax=0.69,
      ymax=1400,
      ]
      \addplot[blue,hist={bins=40}] table [y index = 0,col sep=comma]
      {./data/markdown_bellman_det_vals.csv};
    \end{axis}
  \end{tikzpicture}%
  \begin{tikzpicture}[scale=0.8]
    \begin{axis}[
      xlabel=Profit $P(\alpha^C)$,
      title={Simulations of CEC policy},
      xmin=0.565, xmax=0.69,
      ymax=1400,
      ]
      \addplot[blue,hist={bins=40}] table [y index = 1,col sep=comma]
      {./data/markdown_bellman_det_vals.csv};
      \draw[<-] (axis cs:0.62,220) --
      (axis cs:0.605,600) node[anchor=south] {Fatter lower tail};
    \end{axis}
  \end{tikzpicture}

  \vspace{1em}
  \begin{tikzpicture}[scale=0.8]
    \begin{axis}[
      xlabel=$P(\alpha^B)-P(\alpha^C)$,
      ylabel=Count,
      title={Simulations of Bellman and CEC policies},
      every x tick scale label/.style={font=\boldmath\large,
        at={(xticklabel
           cs:0.98,5pt)},yshift=0em,left,inner sep=0pt,
      }
     ]
      \addplot[red,hist={bins=30}] table [y expr = {(\thisrowno{0}-\thisrowno{1})},col sep=comma]
      {./data/markdown_bellman_det_vals.csv};
      \draw[dashed] (axis cs:0,-30) -- (axis cs:0,6000);
      \draw[<-] (axis cs:0.017,1000) -- (axis cs:0.015,2000)
      node[anchor=south] {$\alpha^B$ best};
      \draw[<-] (axis cs:-0.003,4000) -- (axis cs:0.005,4000) node[anchor=west] {$\alpha^C$ best};
    \end{axis}
  \end{tikzpicture}
  \todo[inline]{Show relative difference
    $(P(\alpha^B)-P(\alpha^C))/P(\alpha^B)$ instead?}
  \caption{This shows the distributions from \num{10000} samples of
    the profits of following the Bellman and CEC policies.
    The sample mean $\mathbb E_W[P(\alpha^B)-P(\alpha^C)]\approx 3.8\times
    10^{-3}$, confirms that Bellman is better on average, as it should be.
    Importantly, however, in more than half the samples the suboptimal policy
    outperforms the Bellman policy.
    The distribution of $P(\alpha^B)-P(\alpha^C)$ appears to be bimodal.
  }\label{fig:bellman_det_vals}
\end{figure}
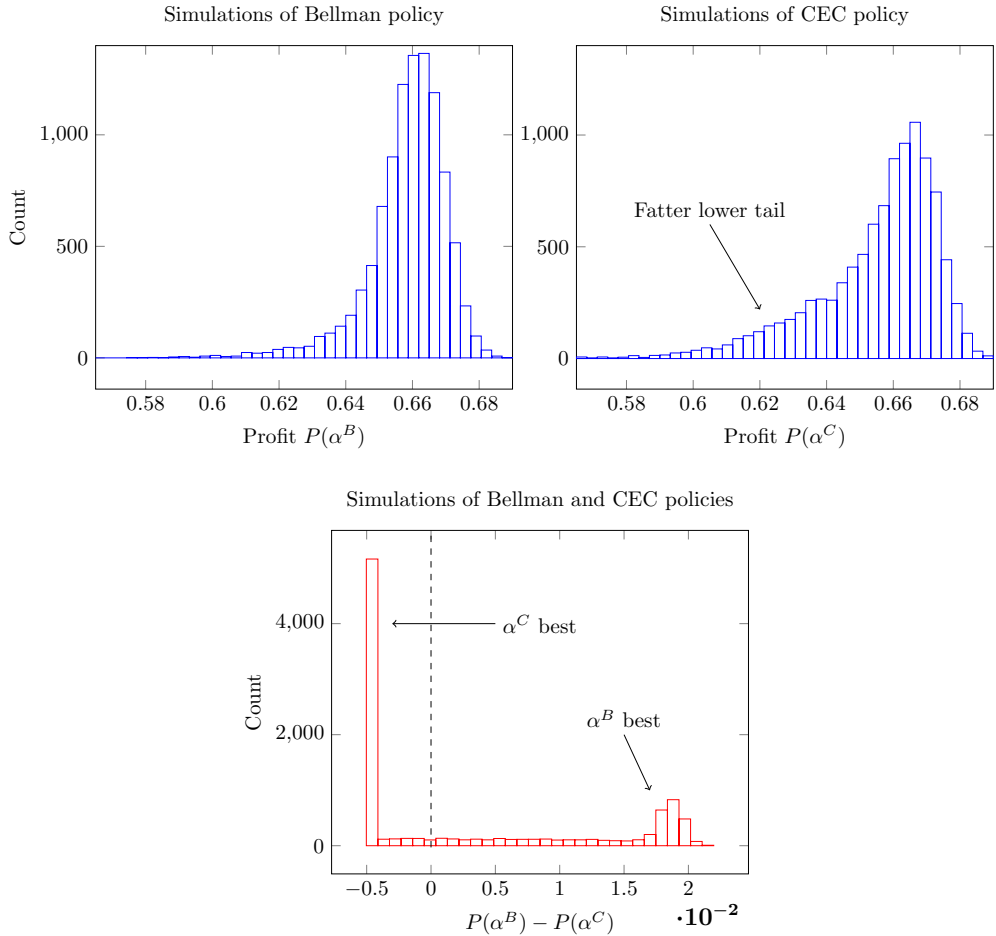

What is more important than how the policy function works, is how it
impacts the goal of the decision-process. Thus, we would like to see
how the two policy processes $\alpha^B$ and $\alpha^C$ perform.
One way to evaluate their performance is to look at
the distribution of the profits $P(\alpha^B)$ and $P(\alpha^C)$.
We remind the reader that
$P(\alpha) =
\sum_{t=0}^{T-1}\alpha_tQ(S_t^\alpha,\alpha_t,W_{t+1})-CS_T^\alpha$.
From the optimality of the Bellman policy function, we must have
that $\mathbb E_W[P(\alpha^B)]\geq \mathbb
E_W[P(\alpha^C)]$. Violations of this result can happen due to
numerical errors in approximating $\alpha^B$ and the expectations.
Marginalising a random variable with the expectation operator, however,
loses a lot of information which can be of interest.
An approximation of the distributions of $P(\alpha^B), P(\alpha^C)$
and their difference
$P(\alpha^B)-P(\alpha^C)$, based on \num{10000} realisations of the
underlying $W$, can be seen in \Cref{fig:bellman_det_vals}.
Indeed in the experiment, the average value of following the
Bellman policy is higher than following the CEC policy.
However, from the bottom figure we see that
in more than half of the cases, the CEC policy outperforms the optimal
policy $\alpha^B$. What we can take from this experiment is that
the CEC policy induces a more risk-seeking pricing strategy than
$\alpha^B$: It results in slightly larger profits for a majority of the
realisations of $W$, but at a cost of taking a more significant
reduction in
profits in the remaining realisations.

\todo[inline]{Add OLFC section here with table similar to \Cref{tbl:paramcomparisons}?}

\subsection{Bellman and CEC parameter comparisons}\label{sec:parameter_comparison}
In the previous sections, our numerical experiments have only shown
results for a fixed combination of the five parameters
termination time $T$, unsold items cost $C$, uncertainty $\gamma$, and
demand function parameters $q_1$ and $q_2$.
In this section we explore the differences between the Bellman policy
and the CEC policy for a larger parameter range.
The formula for $a^C$ in~\eqref{eq:cec_policy} indicates that the
initial price is largely determined by the relationship between
$Tq_1$ and $q_2$, and hence we choose to keep $T=3$ fixed whilst
varying $q_1$ and $q_2$.
We are interested in the difference between the profit following
a Bellman policy $\alpha^B$ and a CEC policy $\alpha^C$.
The policy $\alpha^B$ is computed numerically, and
$\alpha^C$ is obtained using the function $a^C$ from
the formula in \Cref{eq:cec_policy}.

\todo[inline]{Move table to here?}
In particular, the difference between the two policies is measured
using the $L^2$-norm with respect to the probability distribution
induced by the disturbance $(W_1,W_2,\dots,W_T)$, that is
\begin{equation}
  \left\|P(\alpha^B)-P(\alpha^C)\right\|_2
  =\sqrt{\mathbb E_W\left[{( P(\alpha^B)-P(\alpha^C) )}^2 \right]}.
\end{equation}

\begin{figure}[htbp]
  \includegraphics[width=\textwidth]{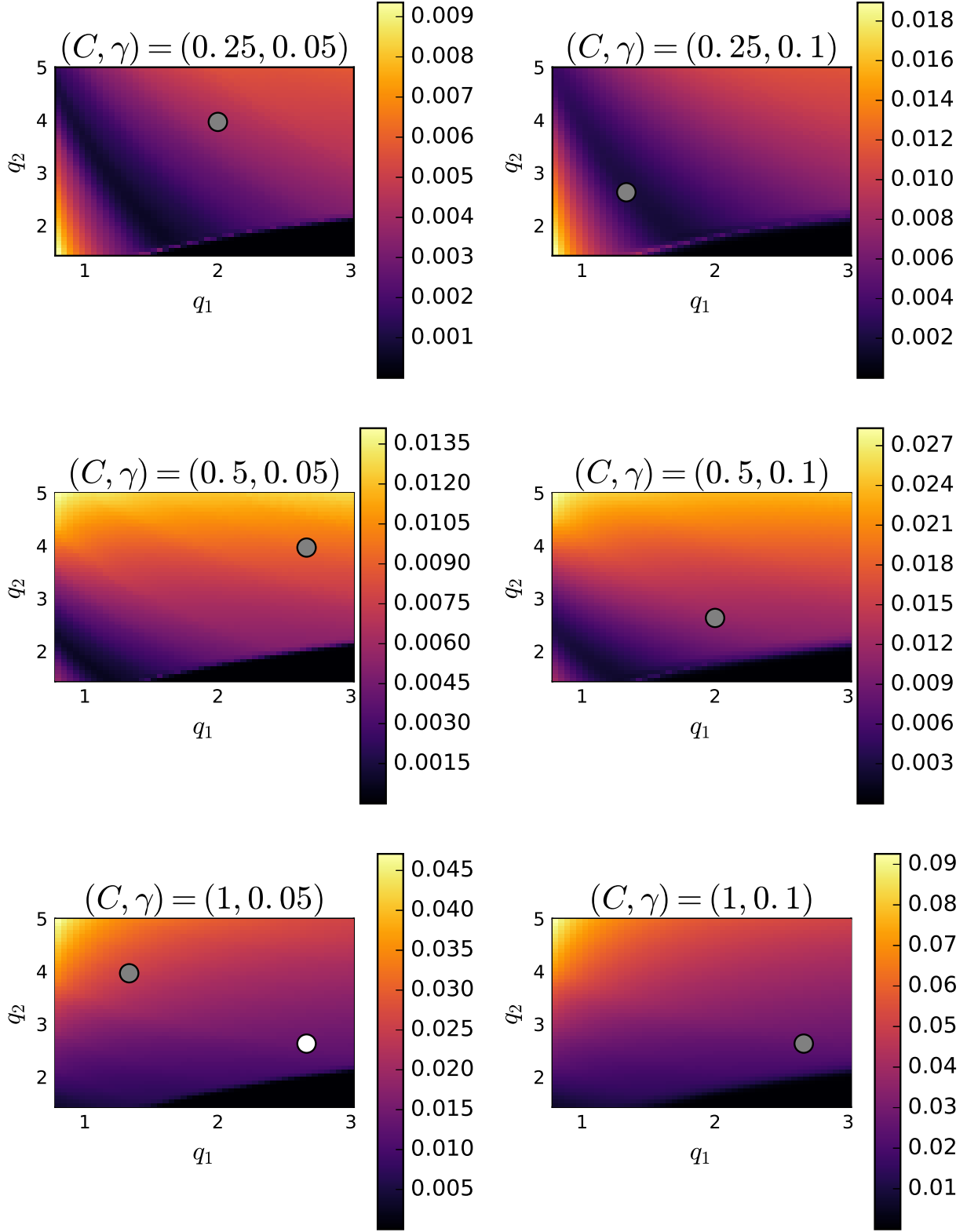}
  \caption{The relative $L^2$ difference
    $\|P(\alpha^B)-P(\alpha^C)\|_2/\|P(\alpha^B)\|_2$ representing
    how different the CEC policy is from the optimal Bellman policy.
    Each plot defines a combination of $(C,\gamma)$, whilst
    the axes vary the parameters $q_1,q_2$ of the demand function
    $q(a)=q_1e^{-q_2a}$.
    The white dot on the bottom, left plot corresponds to the
    parameters chosen for the experiments in this article. The gray
    dots correspond to the parameters in \Cref{tbl:paramcomparisons}.
  }\label{fig:profit_diff_heatmaps}
\end{figure}

To reduce the number of combinations of parameters, we choose only
four combinations of
$(C,\gamma)\in\{(0.5,0.05),(0.5,0.1),(1,0.05),(1,0.1)\}$.
Then,
for each combination of\, $(C,\gamma)$,
we can create a heatmap of the difference between the two policy
functions by varying the parameters $q_1,q_2$.
The relative $L^2$ distance between the CEC and Bellman policy outcomes
is shown in \Cref{fig:profit_diff_heatmaps}. The norms were approximated
with $1000$ samples from $(W_1,W_2,\dots,W_T)$.
Comparing the left and right column, we see that the relative difference doubles as
the standard deviation $\gamma$ is doubled. The cost $C$ plays a large
role, both in terms of the shape of the difference surface and its
magnitude.
The black region at the bottom right of the plots corresponds to popular,
low-elasticity products where both policies suggest to sell at the
maximum price $a=1$.

The parameters used to generate \Cref{fig:bellman_det_vals}
were $C=1$, $\gamma=0.05$, $q_1=e^2/3\approx2.5$, and
$q_2=\nobreak3$. This corresponds to a point
in the region where the relative difference is around $0.016$ --- see the
white dot in the
bottom, left frame in \Cref{fig:profit_diff_heatmaps}.
This difference is in the middle of that seen for all the
combinations of parameters, so the conclusions made in the article
can be considered as relevant for a wider range of systems.
\Cref{tbl:paramcomparisons} provides further data to underscore
the claim that the CEC policy
outperforms the Bellman policy for the
majority of events, at the expense of a stronger underperformance
for the remainder of the events. Notably, the CEC policy is better
more than 50\% of the time for all the parameter
combinations considered in this article.
We can also see from the table that the implicit risk-seeking attitude
of the CEC policy increases with
the relative $L^2$ distance, as its profit distribution widens compared to
the Bellman policy.
In particular, the frequency at which the CEC policy
outperforms the Bellman policy increases, at the expense of a larger
underperformance, or tail loss, in the remaining realisations.
The values in \Cref{tbl:paramcomparisons} were approximated using \num{10000} samples
from $W$, and their corresponding parameter combinations are shown as
grey dots in \Cref{fig:profit_diff_heatmaps}.

\begin{table}[htbp]
  \centering
  \begin{tabular}{llllcccc}
    \toprule
    $C$ & $\gamma$ & $q_1$ & $q_2$ & $\mathcal Q_{0.05}$
    &Median & $\mathcal Q_{0.95}$ &$L^2$\\
    \midrule
    $0.25$ & $0.05$ & $2.0$ & $4.0$  & \textbf{-0.4} & \textbf{-0.3} & \textbf{0.6} & \textbf{0.4}\\
    $0.25$ & $0.1$ & $1.33$ & $2.67$ & \textbf{-0.5} & \textbf{-0.0} & \textbf{0.6} & \textbf{0.3}\\
    $0.5$ & $0.05$ & $2.67$ & $4.0$  & \textbf{-0.6} & \textbf{-0.6} & \textbf{1.9} & \textbf{1.1}\\
    $0.5$ & $0.1$ & $2.0$ & $2.67$   & \textbf{-0.9} & \textbf{-0.6} & \textbf{1.9} & \textbf{1.2}\\
    $1.0$ & $0.05$ & $1.33$ & $4.0$  & \textbf{-1.2} & \textbf{-1.1} & \textbf{5.3} & \textbf{2.7}\\
    $1.0$ & $0.1$ & $2.67$ & $2.67$  & \textbf{-1.5} & \textbf{-1.3} & \textbf{5.8} & \textbf{2.9}\\
        &&&&$\times 10^{-2}$&$\times 10^{-2}$&$\times 10^{-2}$&$\times 10^{-2}$\\
    \bottomrule
  \end{tabular}
  \caption{Statistics comparing profits $P(\alpha^B)$ and $P(\alpha^C)$ for six different
    parameter combinations. The columns $\mathcal Q_s$ represent the
    the $s^{\text{th}}$ quantile of the relative
    difference
    $1-P(\alpha^C)/P(\alpha^B)$. %Median is the same as $\mathcal Q_{0.5}$.
    The values in the column $L^2$ refers to the relative $L^2$
    difference $\|P(\alpha^B)-P(\alpha^C)\|_2/\|P(\alpha^B)\|_2$
    from \Cref{fig:profit_diff_heatmaps}. Each of the parameter
    combinations correspond to a grey dot in \Cref{fig:profit_diff_heatmaps}.
  }\label{tbl:paramcomparisons}
\end{table}

\subsection{Open-Loop Feedback Control policy}
We conclude this article with the example of the Open-Loop
Feedback Control (OLFC) policy~\cite[Ch.~6]{bertsekas2005dynamic},
also known as Stochastic Model
Predictive Control~\cite{farmer2017uncertainty}. This is another suboptimal policy,
similar to the CEC policy, but which better takes into account the uncertainty in the system.
%In fact, one can think of CEC as a zeroth-order approximation to OLFC.\@
The OLFC policy works as follows: At each decision time, a stochastic
optimisation problem for the remaining decision horizon is
solved, based on the most recent quantification of the uncertainty in the
system. Only the decision for the current time step is used, whilst
the subsequent decisions are discarded.
For each decision time $t=0,1,\dots,T-1$, the OLFC policy for the
pricing problem calculates the current price using the following
steps:
\begin{enumerate}
\item Observe the state $s$.
  \begin{samepage}
\item Solve the optimisation problem
  \begin{equation}\label{eq:olfc_optimisation_problem}
    \max_{\mathbf a\in A^{T-t}}\mathbb E_W\left[
      \sum_{\tau=t}^{T-1}\mathbf a_\tau Q(S_\tau^{\mathbf a},\mathbf
      a_\tau,w_{\tau+1})-CS_T^{\mathbf a} \bigm\vert S_t^{\mathbf a}=s \right].
  \end{equation}
\item Set the price corresponding to the element
  $\mathbf{a}_t\in A$ of a maximiser
  to~\eqref{eq:olfc_optimisation_problem}.
  \end{samepage}
\end{enumerate}

For this article, we approximate the expectation operator with
Monte--Carlo using \num{1000} samples, in the same way that the expectation
in the Bellman policy is computed.

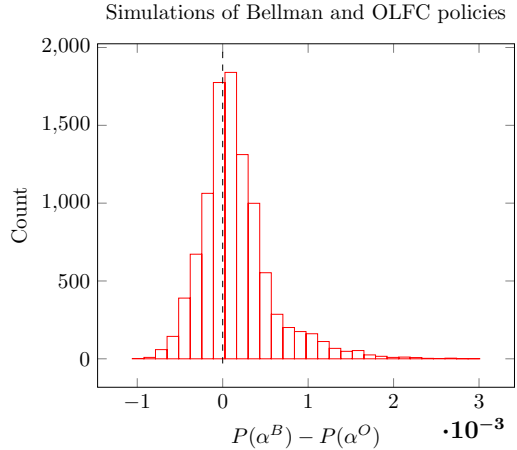
\begin{figure}[htbp]
  \centering
  \begin{tikzpicture}[scale=0.8]
    \begin{axis}[
      xlabel=$P(\alpha^B)-P(\alpha^O)$,
      ylabel=Count,
      title={Simulations of Bellman and OLFC policies},
      every x tick scale label/.style={font=\boldmath\large,
        at={(xticklabel
           cs:0.98,5pt)},yshift=0em,left,inner sep=0pt,
      }
      ]
      \addplot[red,hist={bins=30}] table [y expr = {(\thisrowno{0}-\thisrowno{1})},col sep=comma]
      {./data/markdown_bellman_olfc_vals_1000.csv};
      \draw[dashed] (axis cs:0,-30) -- (axis cs:0,2000);
      % \draw[<-] (axis cs:0.017,1000) -- (axis cs:0.015,2000)
      % node[anchor=south] {$\alpha^B$ best};
      % \draw[<-] (axis cs:-0.003,4000) -- (axis cs:0.005,4000) node[anchor=west] {$\alpha^C$ best};
    \end{axis}
  \end{tikzpicture}
  \todo[inline]{Show relative difference
    $(P(\alpha^B)-P(\alpha^O))/P(\alpha^B)$ instead?}
  \caption{This shows the distributions from \num{10000} samples of
    the profits of following the Bellman and OLFC policies.
    The OLFC policy generates profits an order of magnitude closer to
    the Bellman policy than the CEC policy (\Cref{fig:bellman_det_vals}). This
    distribution appears to be unimodal.
  }\label{fig:bellman_olfc_vals}
  \todo[inline]{Reference}
\end{figure}

As was done for the CEC policy, we compare the performance of following the
OLFC policy $\alpha^O$ with the performance of following the Bellman
policy $\alpha^B$. With the parameters from the example system in
Subsection~\ref{sec:bellman_optimal_control}\ref{sec:bellman_example_markdown},
the empirical distribution of $P(\alpha^B)-P(\alpha^O)$ is shown in
\Cref{fig:bellman_olfc_vals}. The empirical distribution is generated
using \num{10000} samples of $W$.
There are two notable differences in this distribution from the
distribution of $P(\alpha^B)-P(\alpha^C)$ shown in \Cref{fig:bellman_det_vals}.
First, the distribution appears to be unimodal and more concentrated around zero. Second,
the values are an order of magnitude smaller.
This supports the claim that the OLFC policy better approximates the
Bellman policy.

For completeness, \Cref{tbl:paramcomparisons_olfc} includes the statistics of the relative
difference $1-P(\alpha^O)/P(\alpha^B)$ for the same parameters
that were used to compare the CEC and Bellman policies in \Cref{tbl:paramcomparisons}.
The distributions of the relative difference indicate a more symmetric
distribution around the median, with values of an order of magnitude
smaller than in the CEC comparison.
\begin{table}[htbp]
  \centering
  \begin{tabular}{llllcccc}
    \toprule
    $C$ & $\gamma$ & $q_1$ & $q_2$ & $\mathcal Q_{0.05}$
    &Median & $\mathcal Q_{0.95}$ &$L^2$\\
    \midrule
    $0.25$ & $0.05$ & $2.0$ & $4.0$
                                   & \textbf{-0.6} & \textbf{0.1} & \textbf{1.2} & \textbf{0.6} \\
    $0.25$ & $0.1$ & $1.33$ & $2.67$
                                   & \textbf{-2.5} & \textbf{0.4} & \textbf{4.2} & \textbf{2.1} \\
    $0.5$ & $0.05$ & $2.67$ & $4.0$
                                   & \textbf{-0.5} & \textbf{0.1} & \textbf{1.2} & \textbf{0.6} \\
    $0.5$ & $0.1$ & $2.0$ & $2.67$
                                   & \textbf{-2.3} & \textbf{0.3} & \textbf{4.6} & \textbf{2.1} \\
    $1.0$ & $0.05$ & $1.33$ & $4.0$
                                   & \textbf{-0.8} & \textbf{0.0} & \textbf{2.8} & \textbf{1.1} \\
    $1.0$ & $0.1$ & $2.67$ & $2.67$
                                   & \textbf{-2.2} & \textbf{0.4} & \textbf{5.8} & \textbf{2.4} \\
        &&&&$\times 10^{-3}$&$\times 10^{-3}$&$\times 10^{-3}$&$\times 10^{-3}$\\
    \bottomrule
  \end{tabular}
  \caption{Statistics comparing the relative profits of $P(\alpha^O)$
    against $P(\alpha^B)$ for six different
    parameter combinations. Note that the values are an order of magnitude
    smaller than in the comparison between Bellman and CEC shown in \Cref{tbl:paramcomparisons}}\label{tbl:paramcomparisons_olfc}
  \todo[inline]{reference}
\end{table}

Note that the CEC policy is a special case of the OLFC policy, where the
expectation~\eqref{eq:olfc_optimisation_problem} is approximated with~\eqref{eq:cec_optim_problem}, a zeroth-order expansion around
the estimate $w_{t+1},\dots,w_T$.
The OLFC policy is thus much more costly than the CEC policy, but will
also better approximate the optimal Bellman policy. In practice, the
decision maker can balance the cost-versus-optimality by how
accurately they
approximate the expectation~\eqref{eq:olfc_optimisation_problem}.

%%%%%%
% Conclusion
%%%%%%%%%%%%%
\section{Conclusion}\label{sec:conclusion}
In this article we have looked at a mathematical formulation of a
retail pricing problem for profit maximisation, and have investigated the
performance of two algorithms that balance practicality
with degree of suboptimality. The motivation is to better understand
how well the suboptimal policies approximate the returns from the
optimal policy.
Pricing problems are often formulated as an expected value
maximisation, but different algorithms may induce different
distributions of the profits.
Even though the \emph{expected values} of suboptimal policies are not better
than the optimal policy, they may have a higher profit for many
realisations of the underlying probability distribution.
We found in \Cref{sec:suboptimal_approximations} that there are a
large number of
reasonable system parameters for which
the suboptimal Certainty Equivalent Control policy resulted in a higher profit
than the optimal policy in more than half of the realisations.
In the remaining realisations, however, the CEC policy resulted in
much smaller profits. We interpret these results as an indication that
the suboptimal policy is more risk-seeking than the optimal policy, in
a colloquial sense of the term.
The results in this article underscore the importance of looking at
the impact of
different suboptimal algorithms have on the distribution of
an objective, and not only the impact of
marginalised statistics such as the expected value.

The model problem in this article is fairly simple, and
we propose two specific lines for future research that
take this analysis closer to practical models.
First, to investigate multi-product problems where the demand and
availability of one product depends on the other products.
Second, to introduce a state dependence in the uncertainty in the
system, that is, to allow $W_{t+1}$ to depend explicitly on the values
of $S_t$ and $\alpha_t$.
We hope that these two extensions will lead to a better understanding of
whether the effects seen in this article will be stronger or diluted
in real-life problems.

%%%%%%%%%%%%%%%%
% End-stuff
%%%%%%%%%%%%%%%%
\aucontribute{JND and CLF contributed to specifying the problems and
  choosing the types of algorithms to study. ANR devised the method
  for comparing algorithms, performed all the calculations, and drafted
  the first version of the paper. All authors were responsible for the
  final preparation of the manuscript.}
\competing{The authors declare that they have no competing interests.}
\funding{This publication is based on work partially supported by the EPSRC
Centre For Doctoral Training in Industrially Focused Mathematical
Modelling (EP/L015803/1) in collaboration with dunnhumby Limited.}

% Bibliography:
%\clearpage
\bibliographystyle{RS}
\bibliography{references}

\end{document}